\documentclass{article}
\usepackage{amsthm}
\usepackage{amsmath}
\usepackage{amssymb}
\usepackage{amsthm}
\swapnumbers
\theoremstyle{plain}
\newtheorem{theorem}{Theorem}
\newtheorem{proposition}[theorem]{Proposition}
\newtheorem{example}[theorem]{Example}
\newtheorem{remark}[theorem]{Remark}
\title{On the point spectrum of some perturbed \\
differential operators with periodic coefficients }

\author{Igor Cialenco
\thanks{
The author wishes to express his thanks to Prof. P.~Cojuhari for
suggesting the problem and many stimulating conversations, and to
the participants of Analysis Seminar from the Department of
Mathematics, University of Southern California, for their input
and helpful comments.} }

\begin{document}
\maketitle

\begin{abstract}
Finiteness of the point spectrum of linear operators acting in a
Banach space is investigated from  point of view of perturbation
theory. In the first part of the paper we present an abstract
result based on analytical continuation of the resolvent function
through continuous spectrum.  In the second part, the abstract
result is applied to differential operators which  can be
represented as a  differential operator with periodic coefficients
perturbed by an arbitrary subordinated differential operator.

\vskip 0.3 true cm
{\bf Keywords:}
             Spectral theory, perturbation theory, point spectrum,
differential operators, nonselfadjoint operators.
\vskip 0.3 true cm
{\bf  AMS Subject Classifications:}
            47A10, 47A55, 47A75, 47E05, 34L10
\vskip 0.4 true cm

\end{abstract}

\section{Introduction}

One fundamental problem in spectral analysis of linear operators
is finiteness of the point spectrum, which is of great interest
for various  problems from mathematical physics, quantum mechanics
and related topics as well as for spectral theory itself. Unlike
the case of selfadjoint operators, in which  various methods of
investigation were elaborated thanks to the fundamental spectral
theorem, in the case of nonselfadjoint operators this problem is
typically  reduced to the uniqueness theorem of analytic
functions. From this point of view, the problem of finiteness of
the point spectrum (i.e. the set of all eigenvalues including
those contained in the continuous spectrum) of nonselfadjoint
operators has been studied in many papers, where concrete classes
of operators have been considered:
 differential operators (of the second order)  \cite{D, Pavlov},
{S}chr\"odinger operator \cite{Albeverio_Dudkin, Krall},
the Friedrics model  \cite{L1, Stepin},
the finite-difference operators \cite{L2},
perturbed Winer-Hopf integral operators \cite{Kakabadze}.
Some general results have been presented in \cite{Maksudov}, where the nonselfadjoint
operator is considered as a perturbation of a selfadjoint operator acting in some Hilbert space.

In the present paper, using direct methods of perturbation theory, we
propose a generalization of the so-called method of analytic continuation
of the resolvent of the unperturbed operator across its
continuous spectrum for the case of Banach spaces and unbounded operators.
It should be mentioned that the results obtained here are framed in the theory of large perturbations.

In Section 2 the problem of finiteness of the point spectrum
is solved in abstract settings. Conditions on the linear operators $H_0$
and $B$ are given such that the perturbed operator $H_0+B$ has a finite
set of eigenvalues. The
obtained results actually  guarantee  the finiteness of the spectral
singularities which were investigated for the first time by M.~A.~Naimark
\cite{Naimark}, and which are also very closely related to the problem of eigenfunction
expansion and generalized Parseval identity.
The abstract results proposed here can be applied to different classes of
unbounded operators, among which we can mention ordinary differential
operators, integro-differential operators, pseudo-differential operators,
etc. We want to mention that similar problem for bounded (nonselfadjoint)
operators has been studied by
author in \cite{C3, C2,C1}. The abstract results, similar to those presented here,  have been applied to
 finite-difference perturbed operator (of any order),
Wiener-Hopf perturbed operators (abstract, discrete, and integral),
operators generated by Jacobi matrices, etc.

 In Section 3,  we obtain some new results about
finiteness of the point spectrum of the following class of differential operators:
the unperturbed operator is  a
differential operator with periodic coefficients and the perturbation is
a subordinated differential operator. These operators act
in one of the spaces $L_p(\mathbb R)$ or $L_p(\mathbb R_+) \ (1\leq p<+\infty)$.
We want to stress out that while these results are of independent interest,
and could be obtained separately, we derived them from the abstract results presented in Section 2.
Also, we note that the unperturbed operator can be of any order and can be nonselfadjoint itself.
The results agree with know ones. For example, if the unperturbed operator is Hill operator,
then the finiteness of the point spectrum is guaranteed if the potential in the perturbation part
has exponential decay at infinity (see the concluding result of this paper).
%We conclude this paper by considering the perturbed Hill operator.

\section{Abstract Results}
%In this section, we give an abstract scheme on finiteness of the point spectrum
%of a linear operator (generally speaking unbounded) acting in a Banach space .
Let $H_0$ and $B$ be linear operators acting in a Banach space $\mathcal{B}$
such that the following assumptions are fulfilled:
{\it \begin{list}{}{}
\item (i) The spectrum $\sigma(H_0)$ of the operator $H_0$ is
a simple rectifiable curve and the point spectrum of $H_0$ is absent, i.e.
$\sigma_p(H_0)=\emptyset$;
\item (ii) The operator $B$ can be represented in the form $B=RTS$, where
$S$ is an operator acting from $\mathcal{B}$ into $\mathcal{B}_1$
with $\mathrm{Dom}(S)\supset \mathrm{Dom}(H_0)$, the operator $R$
acts from $\mathcal{B}_1$ into $\mathcal{B}$, $T$ acts
in $\mathcal{B}_1$, and $\mathcal{B}_1$ is a Banach space (possibly
different from $\mathcal{B}$).
\end{list} }

We denote by $\mathrm{Dom}(A)$ and $\mathrm{Ran}(A)$ the domain and the range
of the operator $A$, by $\mathbb{B}(\mathcal{B})$ the set of all linear and bounded operators on
$\mathcal{B}$, and by  $\mathbb{B}_\infty(\mathcal{B})$ the class
of all compact operators defined on $\mathcal{B}$. Also, we denote by
$Q(\lambda), \ \lambda\in\mathbb C$,  the operator %of the form
$S(H_0-\lambda I)^{-1}RT$ defined on the set
$\mathcal L_{\lambda} := \{ u\in\mathcal{B}_1 \ :
\ RTu\in\mathrm{Ran}(H_0-\lambda I) \}$.

Under the above assumptions the following statement holds true.
\begin{proposition}\label{prop1}
If $\lambda\in\sigma_p(H_0+B)$, then there exists
$\varphi\in\mathcal L_{\lambda}, \ \varphi\neq 0$, such that
\begin{equation}\label{eq1}
(I+Q(\lambda))\varphi=0 \, .
\end{equation}
\end{proposition}
The proof is based on the following argument. Suppose that $\lambda\in\sigma_p(H_0+B)$.
Then, there exists a vector $u \in\mathcal{B}, \ u\neq 0$, such that
$(H_0-\lambda)u + RTSu=0$. Note that $RTSu\in\mathrm{Ran}(H_0-\lambda I)$, hence
$u+ (H_0-\lambda I)^{-1}RTSu=0$. Consequently, $Su + S(H_0-\lambda I)^{-1}RTSu=0$.
Put $\varphi = Su$ and equality \eqref{eq1} follows. Note that $\varphi\neq 0$, since otherwise
$Su=0$, and then $(H_0 -\lambda I)u=0$, which is a contradiction with initial assumption
$\sigma_p(H_0)=\emptyset$.

It should be mentioned that the operator-valued function
$Q(\lambda)$ plays a key role in perturbation theory and scattering theory.
Proposition similar to Proposition \ref{prop1} show up
in many problems of spectral analysis (see
for instance  \cite{Kuroda2,Kuroda1,Sch} where the selfadjoint case is considered).
Note that \eqref{eq1} does not imply that corresponding $\lambda$ belongs to the point
spectrum  of $H_0$. Actually, $\lambda$ that satisfies \eqref{eq1} is called spectral singularity, and
is related to eigenfunction expansion problem and generalized Parseval identity
(see for instance \cite{Naimark}).

Due to  Proposition \ref{prop1}, to establish that the point spectrum
of the operator $H=H_0+B$ is a finite set, it suffices  to show that the equality
(\ref{eq1}) holds for a finite set of numbers $\lambda\in\mathbb C$, and non-zero
vectors $\varphi\in\mathcal L_{\lambda}$.
Consequently, using the theorem of uniqueness of analytic
operator-valued functions, it is sufficient to establish the analyticity of
the function $Q(\lambda)$ with respect to $\lambda$. Generally speaking,
 $Q(\lambda)$ is analytic only on the resolvent set $\rho(H_0)$ (for instance if
$R,S,T$ are bounded operators), and the  analyticity is lost in the
neighborhood of $\sigma(H_0)$.
In connection with this, we suppose that there exists  an analytic
continuation of the function $Q(\lambda)$ across $\sigma(H_0)$
(of course on  Riemann
surface). Namely, we suppose that the following assumption is satisfied.

Let $\lambda_0\in\sigma(H_0)$, and let $U(\lambda_0)$ be a neighborhood of the point
$\lambda_0$. We denote by $U_0(\lambda_0)$ one of the connected
components of the neighborhood $U(\lambda_0)$ with respect to $\sigma(H_0)$.
In other words, the curve $\sigma(H_0)$ divides the set $U(\lambda_0)$ in several parts, and
by $U_0(\lambda_0)$ we denote the interior of one of these parts.
For example if $H_0$ is a selfadjoint operator, and $\Pi_{\pm}$ denotes
upper/lower half complex plane, then $U_0(\lambda_0) := U(\lambda_0)
\cap \Pi_{+}$ or $U_0(\lambda_0) := U(\lambda_0) \cap \Pi_{-}$.

{\it \begin{list}{}{}
\item
(iii) For every $\lambda_0\in\mathbb C$  there exists a neighborhood
$U(\lambda_0)$ such that for every $U_0(\lambda_0)$ there exists  a
neighborhood $\widehat{U}(\lambda_0)$ (maybe on a Riemann surface) of the point
$\lambda_0$, and  an operator-valued function
$\widehat{Q}(\lambda):\widehat{U}(\lambda_0)\to\mathbb B_{\infty}(\mathcal B_1)$
such that $U_0(\lambda_0)\subset\widehat{U}(\lambda_0), \  \widehat{Q}(\lambda)$
is analytic on $\widehat{U}(\lambda_0)$ and
$\widehat{Q}(\lambda)\supset Q(\lambda), \ \lambda\in
U_0(\lambda_0)$. \end{list} }

\begin{theorem}\label{th} If the operators $H_0$  and $B$ satisfy conditions
(i)-(iii), then the perturbed operator $H=H_0+B$ has a finite set of eigenvalues.
Moreover, the possible eigenvalues have finite multiplicity.
\end{theorem}
\begin{proof} By Assumption (iii), $\widehat{Q}$ is uniquely defined on entire Riemann surface.
Moreover, one can formally write $\widehat{Q}(\lambda)\supset Q(\lambda), \ \lambda\in\mathbb{C}$,
meaning that for every $\lambda\in\mathbb{C}$, there exists $\mu$ on the Riemann surface, such
that $Q(\lambda)=\widehat{Q}(\mu)$.
According to the theorem about holomorphic
operator-valued functions with values in $\mathbb{B}_\infty$ (see, for example, \cite{Kato},
Chapter VII.1.3 or \cite{R-S}, theorem XII.13), the  function
$\widehat{Q}(\cdot)$ has a finite number of zeros in every neighborhood on the Riemann surface.
Consequently $Q(\cdot)$ has a finite number of zeros in $\mathbb{C}$. The possibility of
existence of   sequence  $\lambda_n\in\mathbb{C}, \ \lambda_n\to\infty, \ Q(\lambda_n)=0$, is ruled
out by the analyticity of the resolvent function $(H_0-\lambda I)^{-1}$. Thus,
\eqref{eq1} is satisfied for a finite number of values $\lambda$. Moreover,  for every $\lambda$,
the  subspace $\mathcal{B}_{1}(\lambda)$ generated by the corresponding vectors
$\varphi\in\mathcal{L}_{\lambda}$ that satisfy \eqref{eq1} has
a finite dimension. By Proposition \ref{prop1},
$\sigma_p(H_0+B)$ is finite, and every eigenvalues has finite multiplicity. Theorem is proved.
\end{proof}

In many applications, usually $H_0$ is a selfadjoint operator, so $\sigma(H_0)\subset\mathbb{R}$.
Hence, (i) is satisfied if $H_0$ has no eigenvalues. Condition (ii) is a technical condition,
but strongly related to (iii). In particular, (ii) and (iii) holds true if one may find the operators
$R$ and $S$ such that $R(H_0-\lambda I)^{-1}S$ has analytic continuation
and $T:= S^{-1}BR^{-1}$ is a compact operator. The hardest to check is condition (iii), and
verification depends on the class of  operators to be considered, and the general rule is to have an
explicit or manageable form of the resolvent function of the unperturbed operator.

\section{Application to Differential Operators}
In this section we will present one application of the general results from Section 2. We will consider
some perturbations of differential operators with periodic coefficients of arbitrary order acting
in $L_p(\mathbb{R})$.

Let $H$ be the differential operator of the following form
\begin{equation}\label{eq2}
H\varphi(t)=\sum\limits_{k=0}^{n}h_{k}(t)\frac{d^k\varphi(t)}{dt^k} \, ,
\end{equation}
where $h_k(t)=a_k(t)+q_k(t) \ \ (k=0,1,\dots n; \ t\in\mathbb R \  \mathrm{or}  \
\mathbb R_+), \ \
a_k(t) \ \ (k=0,1,\dots,n;\,a_n\equiv 1)$  are periodic functions (with the same
period $T$) and $q_k(t) \ \ (k=0,1,\dots,n;\,q_n\equiv 0)$  are functions
(generally speaking, complex-valued) vanishing for $t\to\infty.$
Assume that $a_k(t) \ \ (k=0,1,\dots,n)$ are as smooth as required. The
operator $H$ is supposed to act in the spaces $L_p(\mathbb R)$ or $L_p(\mathbb R_+)
\ (1\leq p<+\infty)$. The domain of
the operator $H$ consists of all functions
$\varphi\in L_p(\mathbb R) \ (L_p(\mathbb R_+))$ having
absolutely  continuous derivatives of order $n-1$  on each bounded
interval of the real axis (semiaxis) and derivative  of the $n$-th
order belonging to $L_p(\mathbb R) \ (L_p(\mathbb R_+))$.

To apply the abstract scheme from Section 2, we consider
the operator $H$  as a perturbation of the operator
$$
H_0=\sum\limits_{k=0}^{n}a_{k}(t)\left(\frac{d}{dt}\right)^k
$$
by the differential operator
$$
B=\sum\limits_{k=0}^{n-1}q_{k}(t)\left(\frac{d}{dt}\right)^k.
$$

The spectral properties of the unperturbed operator $H_0$ have been investigated by many
authors (see for instance  \cite{Mc, Rofe} and the references therein).
In \cite{Rofe} the operator $H_0$ is
considered in the space $L_2(\mathbb R)$, while in \cite{Mc} in $L_p(\mathbb R)\
(1\leq p\leq \infty)$. In these papers it is shown  that the spectrum of
the operator $H_0$ is continuous, coincides with the set
 of those $\lambda\in\mathbb C$ for which  the
equation $H_0\varphi=\lambda\varphi$   has  a non-trivial solution
(so-called zones of relative stability), and it is bounded in $\mathbb{C}$.
Moreover, the unperturbed operator $H_0$ has no
eigenvalues and satisfies condition (i) from Section 2. Actually this statement will
also follow from our derivations related to properties (ii) and (iii).
In what follows we suppose that the
operator $H$ is acting in the space $L_p(\mathbb R)$, but all results
(with obvious changes) hold true for $L_p(\mathbb R_+)$.
%It should be noted that the problem of absence of the point spectrum of the
%operator (\ref{eq2}) was studied in \cite{C}.

As we mentioned before, the key point in our abstract scheme is to find an
analytical extensions of function $Q(\lambda)$, for which we need to have at hand a manageable
representation of the resolvent function $(H_0-\lambda I)^{-1}$.
Although the spectrum of the operator  $H_0$ is well-known
(see for instance \cite{Mc, Rofe}),  we will present here  a different method for
describing explicitly the resolvent of $H_0$, suitable for our goal to verify
the abstract conditions (ii) and (iii). The representation relies on
 Floquet-Liapunov theory
about linear differential equations with periodic coefficients
(see for instance \cite{Hartman,Y-S}).

Without loss of generality we can assume that $T=1$.

Let us consider the equation
\begin{equation}\label{eq5}
H_0\varphi=\lambda\varphi \, ,
\end{equation}
where $\lambda$  is a complex number, or in  vector form
\begin{equation}\label{eq6}
\frac{dx}{dt}=A(t,\lambda)\,x \, ,
\end{equation}
where
$$
A(t,\lambda)=\left( \begin{array}{cccccc}
0 & 1 & 0 & \dots & 0  & 0 \\
0 & 0 & 1 & \dots & 0  & 0 \\
\dots & \dots & \dots & \dots & \dots   & \dots \\
0 & 0 & 0 & \dots & 0  & 1  \\
\lambda-a_0& -a_1 & -a_2 & \dots & -a_{n-2}& -a_{n-1}
\end{array}\right) , \ \
x=\left( \begin{array}{c}
\varphi \\
\varphi\prime  \\
\vdots \\
\varphi^{(n-1)} \end{array}
\right).
$$

Denote by $U(t) \ (=U(t,\lambda))$  the matriciant of the equation (\ref{eq6}),
i.e., the matrix which satisfies the following system of differential equations
$$
\frac{dU(t)}{dt}=A(t,\lambda)\,U(t), \ \ U(0)=E_n \, ,
$$
where $E_n$ is $n\times n$ identity  matrix.
The matrix $U(1)$  is called the monodromy matrix of the equation (\ref{eq6}) and
the eigenvalues $\rho_1(\lambda),\dots,\rho_m(\lambda)$ of the matrix $U(1)$ are called
the multiplicators. Also, we will say that
$U(1)$  is the monodromy matrix and
$\rho_1(\lambda),\dots,\rho_m(\lambda)$   are multiplicators
of the operator $H_0-\lambda I$.

Let $\Gamma=\ln U(1)$ be a solutions of the matrix equation
$\exp(\Gamma)=U(1)$. Note that this equation has solutions since
the monodromy matrix $U(1)$ is nonsingular.
Due to Floquet  theory (see for instance \cite{Hartman}), the matrix $U(t)$  has the following representation
\begin{equation}\label{eq9}
U(t)=F(t)\exp(t\Gamma) \, ,
\end{equation}
where $F(t)$  is a nonsingular differentiable matrix of period $T=1$.

Let us describe  explicitly the structure of the matrix $\exp(t\Gamma)$. For this,
we write $\Gamma$ in its Jordan canonical form,
 $\Gamma=GJG^{-1}$,  where
$J=\mathrm{diag}[J_1,\dots,J_m]$ and $J_k, \ k=1,\dots,m$,
is the Jordan canonical  block corresponding to the eigenvalue $\mu_k$.
Hence,
$ \exp(t\Gamma)=G\,\exp(tJ)\,G^{-1}, \ \exp(tJ)=\mathrm{diag}
[\exp(J_1t),\dots,\exp(J_mt)]$,  and
$$
\exp(J_kt)= \exp(t\mu_k)\left(  \begin{array}{cccc}
1 & t  & \dots & \frac{t^{p_k-1}}{(p_k-1)!} \\
0 & 1  & \dots & \frac{t^{p_k-2}}{(p_k-2)!} \\
. & .  & \dots &  \dots \\
0 & 0  & \dots &    1
\end{array}
\right) \, ,
$$
where $p_k, \ k=1,\dots,m$,  is the  dimension of the Jordan block $J_k$.

Since $U(t)$  is the matriciant of the equation (\ref{eq6}), it follows that every
solution of this equation has  the form
\begin{equation}\label{eq12}
x(t)=U(t)\,x_0 \, ,
\end{equation}
where $x_0$  is an arbitrary vector from $\mathbb R^n$.

Thus, from (\ref{eq6})-(\ref{eq12}), we conclude that the components of the vector $x(t)$,
and consequently  the solution of the equation (\ref{eq5}), are linear combinations
of $\exp(\mu_kt) \ \ (k=1,\dots,m)$ with polynomial coefficients.

Remark that $|t^k\exp(t\lambda_j)|\to\infty, \ t\to\infty$,  if
$\mathrm{Re}\lambda_j\not= 0, \  k=0,1,\dots$,  or $\mathrm{Re}\lambda_j=0, \ k=1,2,\dots$.
Also note that $|t^k\exp(t\lambda_j)|=1$  for $k=0, \ \mathrm{Re} \lambda_j=0$.
From this we conclude that the only  solution of equation (\ref{eq5})
belonging to $L_p(\mathbb R) \ (1\leq p<+\infty)$ is the function $\varphi\equiv 0$, which
yields that $\sigma_p(H_0)=\emptyset$.

Since $\rho_j=\exp(\mu_j), \
j=1,\dots,m$, we conclude that $|\rho_j|<1, \  |\rho_j|>1$ or $|\rho_j|=1$,
if and only if $\mathrm{Re}\lambda_j<0,\  \mathrm{Re}\lambda_j>0$ or
$\mathrm{Re}\lambda_j=0$, respectively.

Now we are ready to solve explicitly
equation $H_0u-\lambda u=v$, where $v\in \mathrm{Ran}(H_0-\lambda I)$.
In matrix form this equation becomes
\begin{equation}\label{eq13}
\frac{dx}{dt}=A(t,\lambda)\,x+f \, ,
\end{equation}
where $f=(0,0,\dots,v)^{\bot}, \ \ A(t,\lambda)$ and $x$ are the same as in (\ref{eq6}), and
$\bot$ stands for the transposed vector.
According  to the Floquet representation  of the matriciant (\ref{eq9}), and making
the substitution $x=F(t)\,y$ in (\ref{eq13}), we get
\begin{equation}\label{eq14}
\frac{dy}{dt}=\Gamma\,y+F^{-1}(t)f \, .
\end{equation}
Assume that $|\rho_k|\not= 1, \ k=1,\dots,m$, and suppose that $\rho_k$ are numbered such
that  $|\rho_k|>1$ for $k=1,\dots,l$, and $|\rho_k|<1$ for $k=l+1,\dots,m$.
Denote by $P_1$ the projection in $L_p^n(\mathbb R)$  of the form
$P_1\,y=(0,0,\dots,y_{j+1},\dots,y_n)$, where $y=(y_1,\dots,y_n)
\in L_p^n(\mathbb R), \ j=p_1+\dots+p_l$, and put $P_2:=I-P_1$.

An easy computation shows  that the vector-valued function
\begin{equation}\label{eq15}
y(t)=\int\limits_{-\infty}^{t}\exp(\Gamma(t-s))\,P_1\,F^{-1}(s)\,f(s)ds -
\int\limits^{+\infty}_{t}\exp(\Gamma(t-s))\,P_2\,F^{-1}(s)\,f(s)ds
\end{equation}
is a solution of the equation (\ref{eq14}).

Since $x(t)=F(t)\,y(t)$, it follows that
\begin{eqnarray}\label{eq16}
x(t) & = & F(t)\,\int\limits_{-\infty}^{t}\exp(\Gamma(t-s))\,
P_1\,F^{-1}(s)\,f(s)ds - \nonumber \\
 & - & F(t)\,\int\limits^{+\infty}_{t}\exp(\Gamma(t-s))\,P_2\,F^{-1}(s)\,f(s)ds \, ,
\end{eqnarray}
and taking into account (\ref{eq13}) one obtains
\begin{eqnarray}\label{eq17}
(H_0-\lambda I)^{-1}v(t) & = & \sum\limits_{r=l+1}^{m}\sum\limits_{k=0}^{p_r}
q_{rk}(t)\!\int\limits_{-\infty}^{t}\!\!\!\exp(\mu_r(t-s))(t-s)^kh_{rk}(s)v(s)ds + \nonumber \\
 & + &\sum\limits_{r=1}^{l}\sum\limits_{k=0}^{p_r}q_{rk}(t)\!\int
 \limits^{+\infty}_{t}\!\!\!\exp(\mu_r(t-s))(t-s)^kh_{rk}(s)v(s)ds ,
\end{eqnarray}
where $v\in\mathrm{Ran}(H_0-\lambda I), \ q_{rk}$ and $h_{rk}$ are
some continuous periodic functions.

It is easy to show  that under assumption  $|\rho_k|\not= 1,\ k=1,\dots,m$,
the operator defined in (\ref{eq17}) is bounded, hence $\lambda\in\rho(H_0)$.
Moreover, $\lambda\in\sigma(H_0)$ if there
exists at least one multiplicator which lie on the unit circle $\mathbb T=\{ z\in\mathbb C : |z|=1 \}$.
It should be mentioned, since $a_n(t)\equiv 1$,
it is not possible to have  a  multiplicator $\rho_k$ that belongs to $\mathbb T$ and is independent  of $\lambda$
%$or all the multiplicators $\rho_1,\dots,\rho_m$ are constant and $\rho_k\notin \mathbb T, k=1,\dots,m$
(see for instance \cite{Rofe}).

\begin{remark}\label{remark1}
Summing up, we conclude: the point spectrum of the unperturbed operator
$H_0$ is absent; $\sigma(H_0)$  consists from  the set of all curves
determined by the equation
$\mathrm{det}\left( U(1,\lambda)-\rho I\right)=0, \ |\rho|=1$;
for every regular point $\lambda\in\rho(H_0)$ the resolvent function
$(H_0-\lambda I)^{-1}$ has the form \eqref{eq17}.
\end{remark}

Now we are ready to prove the main result of this section.
In what follows we will preserve
the same notations as we defined above.

\begin{theorem}\label{th1} If the functions $q_k(t), \ k=0,1,\dots,n-1$, are such that
\begin{equation}\label{eq19}
q_k(t)\,\exp(\tau |t|)\in L_{\infty}(\mathbb R)\, ,
\end{equation}
for some $\tau>0$, then the point spectrum of the perturbed operator $H$
is at most  a finite set. Furthermore, the possible eigenvalues have finite
multiplicity.
\end{theorem}
\begin{proof}
 For an arbitrary $\lambda\in\mathbb{C}$, in the space ${L_{p}^{n}}(\mathbb R)$  we consider the operator $H_1$ of the
following form
\begin{equation}\label{eq20}
H_1x(t)=\left( \frac{d}{dt}-A(t,\lambda)\,\right) x(t)+B_1 x(t) \, ,
\end{equation}
where $A(t,\lambda)$ and $x$ are as in (\ref{eq6}), and
$$  B_1= \left(
\begin{array}{cccc}
0 & 0 & \dots & 0 \\
  &   &       &    \\
. . & . . & \dots & . . \\
  &  &       &  \\
q_0(t) & q_1(t) & \dots & q_{n-1}(t)
\end{array}  \right)                        .
$$
Denote by $P_1$ the projection in $L_p^n(\mathbb R)$  of the form
$P_1\,y=(0,0,\dots,y_{l+1},\dots,y_n)$, and put $P_2:=I-P_1$. The index $l$ will be
specified latter on. To satisfy conditions (ii) and (iii) from the
the abstract result, we factorize the operator $B_1$ as follows: $B_1=RTS$,  where
$ R=\exp(-\delta |t|)\cdot P_2$,  $S=\exp(-\delta |t|)\cdot P_1$,
$T=\exp(\delta|t|)\cdot B_1\cdot\exp(\delta|t|)$, and $\delta >0$.
Note that $T\in\mathbb B(L^n_p(\mathbb R))$, and condition (ii) is fulfilled.

By Remark \ref{remark1}, it is  sufficient to show that the
operator-valued function
$$
Q(\lambda)=S\, \left( \frac{d}{dt}-A(t,\lambda)\,\right)^{-1}\!\!R\ T\ ,  \qquad
\lambda\in\sigma(H_0)\, ,
$$
satisfies condition {\it (iii)}.

Let $\lambda_0\in\sigma(H_0)$. Since $U(t)$ is the matriciant,
the matrix  $U(t)=U(t,\lambda)$  is analytic in $\lambda$ (see for instance \cite{Y-S}, p.71, Th.1.3).
Therefore, the function $\mathrm{det}(U(1,\lambda)-\rho)=0$  is also analytic
in $\lambda$. By  implicit  function theorem, there
exists a neighborhood $U(\lambda_0)$ of the point $\lambda_0$, such that
the multiplicators $\rho_1(\lambda),\dots,\rho_m(\lambda)$  are holomorphic
on $U(\lambda_0)$.
Let $U_0(\lambda_0)$ be one of the connected components of the neighborhood
$U(\lambda_0)$, and assume that the multiplicators (including their multiplicity)
are enumerated such that
$$\begin{array}{ccl}
\mathrm{Re} (\mu_k)>0, & \qquad  &  \mathrm{for \ all} \  \lambda\in U_0(\lambda_0);
\,k=1,\dots,l ,  \\
\mathrm{Re} (\mu_k)<0, & \qquad & \mathrm{for \ all} \ \lambda\in U_0(\lambda_0);
\,k=l+1,\dots,m,
\end{array}
$$
where $\exp(\mu_k)=\rho_k \ \ (k=1,\dots,n)$. Using (\ref{eq16}), (\ref{eq17}) and (\ref{eq20}),
we  conclude that the operator $Q(\lambda), \ \lambda\in U_0(\lambda_0)$,
is a linear combination of the following operators
\begin{eqnarray}\label{eq22}
\left(Q_1(\lambda)\varphi\right)(t) & = &
\int\limits^{t}_{-\infty}\exp(\mu_k(\lambda)(t-s))(t-s)^r\exp
(-\tau|t|)\,T\varphi(s)\,ds, \nonumber \\
&  & \qquad \qquad \qquad (k=l+1,\dots,m; \ r\in\mathbb N)  \\
\left(Q_2(\lambda)\varphi\right)(t)& = &
\int\limits_{t}^{+\infty}\exp(\mu_k(\lambda)(t-s))(t-s)^r\exp
(-\tau|s|)\,T\varphi(s)\,ds, \nonumber \\
&  & \qquad  \qquad \qquad (k=1,\dots,l; \ r\in\mathbb N). \nonumber
\end{eqnarray}

We take the neighborhood $\widehat{U}(\lambda_0)\subset U(\lambda_0)$ such that
$\mathrm{Re}(\mu_k(\lambda))-\delta<0, \,  k=1,\dots,l; \ \lambda\in \widehat{U}(\lambda_0)$
and
$\mathrm{Re}(\mu_k(\lambda))+\delta>0, \,  k=l+1,\dots,m;\,\lambda\in \widehat{U}(\lambda_0)$.
For every $\lambda\in\widehat{U}(\lambda_0)$, we define the operator $\widehat{Q}(\lambda)$
by the same formula by which the operator $Q(\lambda)$ is defined on
$U_0(\lambda_0)$ (i.e. integral operators generated by (\ref{eq22})).
Under these conditions, the operator-valued functions
(\ref{eq22}) are analytic on $\widehat U(\lambda_0)$ and take values in $\mathbb{B}_\infty(L_p(\mathbb{R}))$.
Hence,  the same property holds true for the operator-valued function
$\widehat{Q}(\lambda), \, \lambda\in\widehat{U}(\lambda_0)$. By
the definition of $\widehat Q$ we have
$\widehat Q(\lambda)\supset Q(\lambda), \ \lambda\in U_0(\lambda_0)$.
Thus, the condition {\it (iii)} of the abstract scheme is verified,
and Theorem \ref{th1} is proved.
\medskip{}
\end{proof}

\begin{remark}
The initial spectral problem has been reduced to the
corresponding system of  first order differential equations
(\ref{eq6}) and (\ref{eq20}). Moreover, we did not use the particular form and dimension of the matrices
$A(t,\lambda)$ and $B_1(t)$. Actually, the Theorem \ref{th1} holds true
for any matrices $A(t,\lambda)$ and $B_1(t)$, under condition
that $A(t,\lambda)$ is periodic in $t$ and analytic in $\lambda$, and the elements of the matrix $B_1(t)$
are such that $\exp(\tau|t|)|b_{jk}(t)|\in L_\infty(\mathbb R), \
\tau>0;\, j,k=1,\dots,n$. The obtained results are also true if
the operator (\ref{eq2}) is a differential operator with
matrix coefficients, i.e. $a_k(t) \ (k=0,1,\dots,n)$ are periodic matrix-valued
functions of dimension $r\times r, \ \det a_n(t)\neq 0, \ b_k(t), \ k=0,1,\dots,n-1$,
are measurable matrix-valued functions of the same dimension $r\times r$, and the
operator $H$ acts in the space $L_p(\mathbb R,\mathbb C^n) \ (1\leq p<+\infty)$.
In addition, the condition (\ref{eq19}) should be replaced by the following one:
$\exp(\tau|t|)|b_k|\in L_\infty(\mathbb R)$,
for some $\tau>0$, where $|\cdot|$ is the operator matrix norm in $\mathbb C^n$.
\end{remark}

\begin{remark}
It should be mentioned that  similar results can be
formulated for more general classes  of operators. Namely, instead of periodicity of
the unperturbed operator it is sufficient to suppose that the system of
differential equations generated by the unperturbed operator is a reducible one.
\end{remark}
\begin{example}
 As a concrete application of the previous results, consider the
differential operator of the form
\begin{equation}\label{eq23}
H=-\frac{d^2}{dt^2}+q_1(t)\frac{d}{dt}+p(t)+q_2(t) \, ,
\end{equation}
where $p(t+1)=p(t), \ q_k, \ k=1,2$, are measurable, complex-valued functions,
and the operator $H$ is acting in $L_2(\mathbb R)$. The unperturbed operator is
well-studied Hill operator (see for instance \cite{Glazman}, p.281)
$$
H_0u=-\frac{d^2u}{dt^2}+p(t)u \, .
$$
Hence, for the perturbed Hill operator (\ref{eq23}) we can formulate the following
result.

\begin{theorem}
 If $\exp(\tau|t|)q_k(t)\in L_\infty(\mathbb R) \ (k=1,2;\,\tau>0)$,
then the perturbed Hill operator (\ref{eq23})
has a finite set of eigenvalues, each of them of finite multiplicity.
\end{theorem}
\end{example}

\def\cprime{$'$} \def\cprime{$'$} \def\cprime{$'$}
\providecommand{\bysame}{\leavevmode\hbox
to3em{\hrulefill}\thinspace}
\providecommand{\MR}{\relax\ifhmode\unskip\space\fi MR }
\providecommand{\MRhref}[2]{%
  \href{http://www.ams.org/mathscinet-getitem?mr=#1}{#2}
} \providecommand{\href}[2]{#2}

\medskip
Igor Cialenco \hfill \break
Illinois Institute of Technology \hfill  \break
Department of Applied Mathematics \hfill  \break
10 W 32nd Str, E1 234C   \hfill  \break
Chicago, IL 60659
  \hfill\break
e-mail:\ igor@math.iit.edu    \hfill

\end{document}